\documentclass[11pt]{amsart}
\usepackage{amsmath}
\usepackage{cases}
\usepackage{mathrsfs}
\usepackage{bbm}
\usepackage{amssymb}
\usepackage{txfonts}
\usepackage{amscd}
\usepackage{amsfonts,latexsym,amsmath,amsthm,amsxtra,mathdots,amssymb,latexsym,mathabx}
\usepackage[all,cmtip]{xy}
\RequirePackage{amsmath} \RequirePackage{amssymb}
\usepackage{color}
\usepackage{colordvi}
\usepackage{multicol}
\usepackage{hyperref}

\usepackage[margin=1in]{geometry}
\usepackage{xcolor}
\hypersetup{
    colorlinks,
    linkcolor={red!50!black},
    citecolor={blue!50!black},
    urlcolor={blue!80!black}
}

\def\mod{\mathrm{mod}\ }

\def\sumx{\sideset{}{^\star}\sum}

    \newcommand{\BC}{{\mathbb {C}}}

     \newcommand{\BR}{{\mathbb {R}}}

     \newcommand{\BZ}{{\mathbb {Z}}}

\newcommand{\Z}{\mbox{$\mathbb Z$}}	
\newcommand{\R}{\mbox{$\mathbb R$}}	

\def\-{^{-1}}

\def\-{^{-1}}

\newcommand{\delete}[1]{}

    \theoremstyle{plain}

    \newtheorem{thm}{Theorem}[section] 
    \newtheorem{lem}[thm]{Lemma}  \newtheorem{prop}[thm]{Proposition}

    \newtheorem {rem}[thm]{Remark}

    \numberwithin{equation}{section}

\begin{document}

\title{$t-$ASPECT SUBCONVEXITY FOR $GL(2)-$L FUNCTIONS}
\maketitle
\begin{center}
\author{Keshav Aggarwal}
\end{center}
\begin{abstract}
Let $f$ be a holomorphic cusp form for $SL_2(\Z)$ of weight $k>1$. In these notes, we follow Munshi \cite{munshi} to prove the Burgess bound
$$ L(1/2+it,f) \ll_{f,\varepsilon} (1 + |t|)^{1/2 - 1/8 +\varepsilon}. $$
\end{abstract}

\section{Introduction}

Let $f$ be a holomorphic cusp form for $SL_2(\Z)$ of weight $k>1$. The $L$-series is given by,
$$ L(s,f) = \sum_{n\geq1} \lambda_f(n) n^{-s} \quad \textit{ for } Re(s)>1. $$ 

This extends to an entire function on the whole complex plane $\BC$. The convexity principle gives the bound $L(1/2+it,f)\ll_f(1+|t|)^{1/2}$, known as the convexity bound. The purpose of this paper is to prove the following bound.
\begin{thm}\label{mainthm}
Let $f$ be a holomorphic cusp form for $SL_2(\BZ)$. Then we have,
$$ L(1/2+it,f) \ll_{f,\varepsilon} (1 + |t|)^{1/2 - 1/8 +\varepsilon}. $$
\end{thm}
The first such bound was obtained by Good \cite{good82}. The result was extended to Maass cusp forms by Jutila \cite{ju97}. $t$-aspect subconvexity for higher $GL(n)$ is largely unknown. Subconvex bounds for $GL(1)$ and $GL(2)$, uniformly in all aspects is known by the works of Michel-Venkatesh \cite{mive10}. $t$-aspect subconvexity for self dual Hecke-Maass forms for $GL(3)$ was first established by Li \cite{li11}. Munshi \cite{munshi} used a different method (that we follow and execute) to extend the result to all Hecke-Maass cusp forms. Recently, Singh \cite{singh17} did similar calculations for $t$-aspect subconvexity for $GL(2)$ $L$-functions of holomorphic and Hecke-Maass cusp forms and claims to get the Weyl bound.

We have followed the ideas of Munshi \cite{munshi} and use a modification of the circle method. In the present situation, Kloosterman's version of the circle method works best. Let,
\[
\delta(n)=
\begin{cases}
&1 \quad \textit{ if }n=0,\\
&0 \quad \textit{ if }n\neq0.
\end{cases}
\]
Then for any real number $Q>0$, we have,
\begin{equation}\label{delta}
\delta(n)=2Re\int_0^1\underset{1\leq q\leq Q<a\leq Q}{ \sum \sideset{}{^*}\sum} \frac{1}{aq}e\left(\frac{n\overline{a}}{q}-\frac{nx}{aq}\right)dx
\end{equation}
for $n\in\BZ$. Here $e(.)=e^{2\pi i.}$ and the $*$ on the inner sum means that $(a,q)=1$. $\overline{a}$ is the multiplicative inverse of $a\mod q$. There are well understood drawbacks of this circle method. It will turn out that this circle method in itself will not be sufficient, and we will have a apply a `conductor lowering trick' as used by Munshi in his various works \cite{munshi, mu15}.

Suppose $t>2$. The approximate functional equation gives
$$ L(1/2+it,f)\ll t^\epsilon \sup_{N\leq t^{1+\epsilon}}\frac{|S(N)|}{N^{1/2}} + t^{-2015} $$
where
$$ S(N) := \sum_{n\geq1} \lambda(n) n^{-it} V\left(\frac{n}{N}\right). $$

Let $V$ be a smooth function supported on $[1,2]$ satisfying $V^{(j)}\ll_j1$. We further normalize $V$ so that $\int_\BR V(x)dx=1$. We will apply (\ref{delta}) directly to $S(N)$ with a conductor lowering integral to separate the oscillations of $\lambda(n)$ and $n^{-it}$.
\begin{equation}
S(N)=\frac{1}{K}\int V\left(\frac{v}{K}\right)\sum_{n\geq1}\sum_{m\geq1}\lambda(n)m^{-it}\left(\frac{n}{m}\right)^{iv}V\left(\frac{n}{N}\right)U\left(\frac{m}{N}\right)\delta(n-m)dv.
\end{equation}
where $t^\varepsilon<K<t$ is a parameter that will be optimized later. $U$ is a smooth function which is supported on $[1/2,5/2]$, with $U(x)=1$ on $[1,2]$ and satisfies $U^{(j)}\ll_j1$. The extra integral introduced is
\begin{equation*}
\frac{1}{K}\int \left(\frac{n}{m}\right)^{iv}V\left(\frac{v}{K}\right)dv.
\end{equation*}
For $n, m\in[N,2N]$, integration by parts shows that the above integral is small if $|n-m|\gg Nt^{\varepsilon}/K$. This is the crucial `trick' in the paper. As Munshi points out in the $SL_3(\BZ)$ case \cite{munshi}, introduction of this parameter $K$ will seem to hurt us until the very last step, which we will justify in the proof sketch.

We can therefore write $S(N) = S^+(N) + S^-(N)$ where
\begin{equation}\label{start}
\begin{split}
S^\pm(N) = & \frac{1}{K} \int_0^1\int_{\R} V\left(\frac{v}{K}\right) \underset{1\leq q\leq Q<a\leq Q}{ \sum \sideset{}{^*}\sum} \frac{1}{aq} \sum_{n,m\geq 1} \lambda(n) n^{iv} m^{-i(t+v)}\\ & e\left(\pm\frac{(n-m)\bar{a}}{q} \mp \frac{(n-m)x}{aq}\right) V\left(\frac{n}{N}\right) U\left(\frac{m}{N}\right) dv dx.
\end{split}
\end{equation}
The analysis and bounds for $S^+(N)$ and $S^-(N)$ are similar. We therefore analyze only $S^+(N)$. We will justify later in Remark \ref{conductorlowering} that the natural choice for $Q$ is $Q=(N/K)^{1/2}$ (and thus the lowering of conductor by $K^{1/2}$).

We will take 
\begin{equation}
t^{3/4} \ll N<t^{1+\varepsilon} \textit{ and } N^{1/2}\leq K\ll N^{1-\varepsilon}
\end{equation}
In this range, we will establish the following bound.

\begin{prop}\label{mainprop}
For $t^{3/4} \ll N<t^{1+\epsilon}$, we have
\begin{equation}
\frac{S^+(N)}{N^{1/2}}\ll t^{1/2+\varepsilon}\left(\frac{K^{1/2}}{N^{1/2}}+\frac{1}{K^{1/4}}\right).
\end{equation}
\end{prop}
Same bound holds for $S^-(N)$, and consequently for $S(N)$. The optimal choice of $K$ is therefore $K=N^{2/3}$. With this choice of $K$, $S(N)/N^{1/2}\ll t^{1/2}/N^{1/6}$. For $N\ll t^{3/4}$, the trivial bound $S(N)\ll Nt^{\varepsilon}$ is sufficient. This follows by applying Cauchy's inequality to the $n$-sum followed by Lemma \ref{ramonavg} (Ramanujan bound on average). Theorem \ref{mainthm} then follows from Lemma \ref{ramonavg} and Proposition \ref{mainprop}.

\subsection{Proof Sketch}
We briefly explain the steps of the proof and provide heuristics in this subsection. Temporarily assume Ramanujan conjecture $\lambda(n)\ll n^{\varepsilon}$. This is not a serious assumption, since at any step we can apply Cauchy inequality and use Lemma \ref{ramonavg}. The circle method is used to separate the sums on $n$ and $m$, and we arrive at (\ref{start}). Trivial estimate gives $S(N)\ll N^{2+\varepsilon}$. For simplicity, let $N\asymp t$ and $q\asymp Q$. So we are required to save $N$ and a little more in a sum of the form
\begin{equation*}
\int_K^{2K}\sum_{q\asymp Q}\quad\sumx_{Q<a\leq Q+q}\sum_{n\asymp N}\lambda(n)n^{iv}e\left(\frac{n\overline{a}}{q}-\frac{nx}{aq}\right)\sum_{m\asymp N}m^{-i(t+v)}e\left(\frac{-m\overline{a}}{q}+\frac{mx}{aq}\right)dv.
\end{equation*}
The sum over $m$ has `conductor' $Qt\asymp N^{1/2}t/K^{1/2}$. Roughly speaking, the conductor takes into account the arithmetic modulus $q$, with the size $=(t+v)$ of oscillation of the analytic weight. If we assume $K\ll t^{1-\varepsilon}$, then the size of the oscillation is $t$, so the extra oscillation of $m^{-iv}$ does not hurt us here. Poisson summation changes the length of summation to $Qt/N\asymp Q$, and contributes a factor of $N$ along with a congruence condition mod $q$ and an oscillatory integral. The oscillatory integral saves us $t^{1/2}$. In all, we will save $N/t^{1/2}$ in this step. So far the saving is independent of $K$. Next step is to apply Voronoi summation to the $n$-sum. We need to save $t^{1/2}$ in a sum of the form
\begin{equation*}
\int_K^{2K}\sum_{q\asymp Q}\underset{\begin{subarray}{c}(m,q)=1\\|m|\ll Qt/N\end{subarray}}{\sum}\left(\frac{(t+v)aq}{(x-ma)}\right)^{-i(t+v)}\sum_{n\asymp N}\lambda(n)e\left(\frac{nm}{q}\right)n^{iv}e\left(-\frac{nx}{aq}\right)dv,
\end{equation*}
where $a$ is the unique multiplicative inverse of $m\mod q$ in the range $(Q,q+Q]$. Since the $n$-sum involves $GL(2)$ Fourier coefficients, the `conductor' for the $n$-sum would be $(QK)^2$. The new length of sum would be $(QK)^2/N\asymp K$. Voronoi summation would contribute a factor of $N/q$, a dual additive twist and an oscillatory weight function. The oscillation in the weight function would save us $K^{1/2}$. In all, we will save $Q/K^{1/2}=N^{1/2}/K$. If $K$ is large, we are actually making it worse. We are therefore left to save $t^{1/2}K/N^{1/2}$ in $S(N)$. Using stationary phase analysis, we will be able to save $K^{1/2}$ in the integral over $v$. At this point, $K$ seems to be hurting more than helping. The final step is to get rid of the $GL(2)$ oscillations using Cauchy inequality and then change the structure using Poisson summation formula. After Cauchy, the sum roughly looks like,
\begin{equation*}
\bigg[\sum_{n\ll K}\bigg|\underset{\begin{subarray}{c}q\asymp Q\\ |m|\asymp Qt/N\\(m,q)=1\end{subarray}}{\sum}e\left(-\frac{nm}{q}\right)\int_{-K}^Kn^{-i\tau}g(q,m,\tau)d\tau\bigg|^2\bigg]^{1/2}.
\end{equation*}
where $g(q,m,\tau)$ is an oscillatory weight function of size $O(1)$. The next steps would be to open the absolute value squared and, apply Poisson to the $n$-sum and analyze the $\tau$-integral. The $\tau$-integral gives us a saving of $K^{1/2}$. After Cauchy and Poisson summation, we will save $N^{1/2}/K^{1/2}$ in the diagonal term and $K^{1/4}$ in the off-diagonal term. Saving over convexity bound in the diagonal terms is $N^{1/2}/K^{1/2}$. Saving over convexity from the off-diagonal terms is $K^{1/4}$. We will therefore get maximum saving when $N^{1/2}/K^{1/2}=K^{1/4}$, that is $K=N^{2/3}$. That gives us a saving of $N^{1/6}$ over the convexity bound of $t^{1/2+\varepsilon}$. Matching this with the trivial bound $N^{1/2}$ for $N\ll t^{3/4}$ gives us the Burgess bound.

\section{$GL(2)$ Voronoi formula and Stationary phase method}

\subsection{Voronoi summation formula for $SL_2(\BZ)$} Suppose $f$ is a holomorphic cusp form for $SL_2(\BZ)$ which is an eigenfunction for all Hecke operators with $n^{th}$ Fourier coefficient $\lambda(n)$, normalized so that $\lambda(1)=1$. In this subsection, we will mention two important results- a summation formula for Fourier coefficients twisted by an additive character, and a bound on the average size of these Fourier coefficients, both of which will play a crucial role in our analysis.

Let $F$ be a smooth function compactly supported on $(0,\infty)$, and let $\tilde{F}(s)=\int_0^\infty g(x)x^{s-1}dx$ be its Mellin transform. An application of the functional equation of $L(s,f)$, followed by unwinding the integral and shifting the contour gives the Voronoi summation formula \cite{kmv}.

\begin{lem}\label{voronoi}
\begin{equation}
\sum_{n\geq1}\lambda(n)e\left(n\frac{a}{q}\right)F(n)=\frac{1}{q}\sum_{n\geq1}\lambda(n)e\left(-n\frac{\overline{a}}{q}\right)\int_0^\infty F(x)\left[2\pi i^kJ_{k-1}\left(\frac{4\pi\sqrt{nx}}{q}\right)\right]dx.
\end{equation}
For our calculations, we take a step back and use the following representation of $J_{k-1}$ as an inverse Mellin transform,
\begin{equation}
J_{k-1}(x)=\frac{1}{2}\frac{1}{2\pi i}\int_{(\sigma)}\left(\frac{x}{2}\right)^{-s}\frac{\Gamma(s/2+(k-1)/2)}{\Gamma(1-s/2+(k-1)/2)} \quad \textit{ for }0<\sigma<1.
\end{equation}
\end{lem}

We would need to study the oscillation of the gamma factors more closely. Recall the Stirling's formula,
\begin{equation*}
\Gamma(\sigma+i\tau) = \sqrt{2\pi}(i\tau)^{\sigma-1/2}e^{-\pi|\tau|/2}\left(\frac{|\tau|}{e}\right)^{i\tau}\left\lbrace 1+ O\left(\frac{1}{|\tau|}\right) \right\rbrace
\end{equation*}
as $|\tau|\rightarrow\infty$. Letting $\gamma(s)=\frac{\Gamma(s/2+(k-1)/2)}{\Gamma(1-s/2+(k-1)/2)}$, we get
\begin{equation}
\gamma(1+i\tau)=\left(\frac{|\tau|}{4e\pi}\right)^{i\tau}\Phi(\tau), \quad \textit{ where }\quad\Phi'(\tau)\ll\frac{1}{|\tau|}.
\end{equation}

We would also need the following bound, which gives Ramanujan conjecture on average. It follows from standard properties of Rankin-Selberg $L$-functions and is well known.
\begin{lem}\label{ramonavg}
We have,
\begin{equation*}
\sum_{n\leq x}|\lambda(n)|^2\ll_{f,\varepsilon} x^{1+\varepsilon}.
\end{equation*}
\end{lem}

\subsection{Stationary phase method} We will need to estimate integrals of the type
\begin{equation}\label{spi}
\mathcal{I}=\int_a^bg(x)e(f(x))dx.
\end{equation}
Let supp$(g)\subset[a,b]$ and $g^{(j)}(x)\ll_{j,a,b}1$. Further suppose there is a $B>0$ such that for $x\in[a,b]$, $|f'(x)|\gg B$ and $f^{(j)}(x)\ll B^{1+\varepsilon}$ when $j>1$. Integration by parts $j$-times gives $\mathcal{I}\ll_{j,a,b,\varepsilon}B^{-j+\varepsilon}$.

In case $f'(x)=0$ at a unique point $x=x_0\in[a,b]$, there is an asymptotic expansion of the integral around $x_0$. $x_0$ is called the stationary phase. A sharp version useful for us can be found in \cite{bky, huxley}.

\begin{lem}\label{stationary}
Suppose $f$ and $g$ are smooth real valued functions satisfying
\begin{equation} \label{spb}
f^{(i)}(x) \ll \Theta_f / \Omega_f^i, \quad g^{(j)}(x) \ll 1/\Omega_g^j
\end{equation}
for $i=2,3$ and $j=0,1,2$. Suppose $g(a)=g(b)=0$. Define
\begin{equation*}
\mathcal{I} = \int_a^b g(x) e(f(x)) dx.
\end{equation*}
\begin{enumerate}
\item Suppose $f'$ and $f''$ do not vanish in $[a,b]$. Let $\Lambda = \min_{[a,b]}|f'(x)|$. Then we have
\begin{equation} \label{SPB}
\mathcal{I} \ll \frac{\Theta_f}{\Omega_f^2\Lambda^3}\left(1 + \frac{\Omega_f}{\Omega_g} + \frac{\Omega_f^2}{\Omega_g^2}\frac{\Lambda}{\Theta_f/\Omega_f}\right).
\end{equation}
\item Suppose $f'$ changes sign from negative to positive at the unique point $x_0\in (a,b)$. Let $\kappa = \min\{b-x_0, x_0-a\}$. Further suppose that (\ref{spb}) holds for $i=4$ and 
\begin{equation}\label{spb1}
f^{(2)}(x) \gg \Theta_f/\Omega_f^2
\end{equation}
holds. Then
\begin{equation}\label{SP}
\mathcal{I} = \frac{g(x_0)e(f(x_0) + 1/8)}{\sqrt{f''(x_0)}} + O\left(\frac{\Omega_f^4}{\Theta_f^2\kappa^3} + \frac{\Omega_f}{\Theta_f^{3/2}} + \frac{\Omega_f^3}{\Theta_f^{3/2}\Omega_g^2}\right).
\end{equation}
\end{enumerate}
\end{lem}

We will also need a second derivative bound for integrals in two variables. Let
\begin{equation}
\mathcal{I}_{(2)}=\int_a^b\int_c^d g(x,y)e(f(x,y))dydx.
\end{equation}
with $f$ and $g$ smooth real valued functions. Let supp$(g)\subset(a,b)\times(c,d)$. Let $r_1, r_2$ be such that inside the support of the integral,
\begin{equation}\label{2nd.der.bound}
f^{(2,0)}(x,y)\gg r_1^2, \quad f^{(0,2)}(x,y)\gg r_2^2, \quad f^{(2,0)}(x,y)f^{(0,2)}(x,y)-\left[f^{(1,1)}(x,y)\right]^2\gg r_1^2r_2^2,
\end{equation}
where $f^{(i,j)}(x,y)=\frac{\partial^{i+j}}{\partial x^i\partial y^j}f(x,y)$. Then we have (see \cite{srini65}),
\begin{equation*}
\mathcal{I}_{(2)}\ll\frac{1}{r_1r_2}.
\end{equation*}
Define the total variance of $g$ by
\begin{equation*}
\textit{var}(g):=\int_a^b\int_c^d\left|\frac{\partial^2}{\partial x\partial y}g(x,y)\right|dydx.
\end{equation*}
Integration by parts along with the above bound gives us the following.
\begin{lem}\label{2nd.der.lem}
Suppose $f, g, r_1, r_2$ are as above and satisfy condition (\ref{2nd.der.bound}). Then we have
\begin{equation*}
\mathcal{I}_{(2)}\ll\frac{\textit{var}(g)}{r_1r_2}.
\end{equation*}
\end{lem}

\subsection{An integral of interest}
Following Munshi \cite{munshi}, let $W$ be a smooth real valued function with supp$(W)\subset[a,b]\subset(0,\infty)$ and $W^{(j)}(x)\ll_{a,b,j}1$. Define
\begin{equation}
W^{\dagger}(r,s)\int_0^\infty W(x)e(-rx)x^{s-1}dx
\end{equation}
where $r\in\BR$ and $s=\sigma+i\beta\in\BC$. This integral is of the form (\ref{spi}) with
\begin{equation*}
g(x)=W(x)x^{\sigma-1} \quad \textit{ and } \quad f(x)=-rx+\frac{1}{2\pi}\beta\log x.
\end{equation*}
Then,
\begin{equation*}
f'(x)=-r+\frac{1}{2\pi}\frac{\beta}{x}\quad\textit{ and }\quad f^{(j)}(x)=(-1)^j(j-1)!\frac{1}{2\pi}\frac{\beta}{x^j}
\end{equation*}
for $j\geq2$. The unique stationary phase occurs at $x_0=\beta/2\pi r$. Note that we can write
\begin{equation}
f'(x)=\frac{\beta}{2\pi}\left(\frac{1}{x}-\frac{1}{x_0}\right)=r\left(\frac{x_0}{x}-1\right).
\end{equation}
Applying Lemma \ref{stationary} appropriately to $W^{\dagger}(r,s)$, we get the following.
\begin{lem}\label{dagger}
Let $W$ be a smooth real valued function with supp$(W)\subset[a,b]\subset(0,\infty)$ and $W^{(j)}(x)\ll_{a,b,j}1$. Let $r\in\BR$ and $s=\sigma+i\beta\in\BC$. We have
\begin{equation}
W^{\dagger}(r,s)=\frac{\sqrt{2\pi}e(1/8)}{\sqrt{-\beta}}W\left(\frac{\beta}{2\pi r}\right)\left(\frac{\beta}{2\pi r}\right)^\sigma\left(\frac{\beta}{2\pi er}\right)^{i\beta}+O_{a,b,\sigma}\left(\min\{|\beta|^{-3/2},|r|^{-3/2}\}\right).
\end{equation}
We also have
\begin{equation}
W^{\dagger}(r,s)=O_{a,b,j,\sigma}\left(\min\left\lbrace\left(\frac{1+|\beta|}{|r|}\right)^j,\left(\frac{1+|r|}{|\beta|}\right)^j\right\rbrace\right).
\end{equation}
\end{lem}

\section{Application of dual summation formulas}

\subsection{Poisson summation to the $m$-sum}
The $m$-sum is given by
\begin{equation*}
\sum_{m\geq1}m^{-i(t+v)}e\left(\frac{-m\bar{a}}{q} + \frac{mx}{aq} \right)U\left(\frac{m}{N}\right) dv dx.
\end{equation*}
Breaking the $m$-sum into congruence classes modulo $q$, we get
\begin{equation*}
\sum_{\alpha\mod q}e\left(\frac{-\alpha\overline{a}}{q}\right)\sum_{m\in\BZ}(\alpha+mq)^{-i(t+v)}e\left(\frac{(\alpha+mq)x}{aq}\right)U\left(\frac{\alpha+mq}{N}\right).
\end{equation*}
Poisson to the $m$-sum gets us
\begin{equation*}
\sum_{\alpha\mod q}e\left(\frac{-\alpha\overline{a}}{q}\right)\sum_{m\in\BZ}\int_{\BR}(\alpha+yq)^{-i(t+v)}e\left(\frac{(\alpha+yq)x}{aq}\right)U\left(\frac{\alpha+yq}{N}\right)e(-my)dy.
\end{equation*}
Making the change of variables $(\alpha+yq)\mapsto u$ and executing the complete character sum $\mod q$, we arrive at
\begin{equation}
N^{1-i(t+v)}\underset{\begin{subarray}{c}m\in\BZ\\m\equiv\overline{a}\mod q\end{subarray}}{\sum}\int_{\BR}U(u)u^{-i(t+v)}e\left(\frac{N(x-ma)}{aq}u\right)du.
\end{equation}
The above integral equals
\begin{equation}
U^\dagger\left(\frac{N(ma-x)}{aq}, 1-i(t+v)\right).
\end{equation}
Everything together,
\begin{equation}
\begin{split}
S^+(N) = &\frac{1}{K} \int_0^1\int_{\R} V\left(\frac{v}{K}\right) \underset{1\leq q\leq Q<a\leq Q}{ \sum \sideset{}{^*}\sum} \frac{1}{aq} \sum_{n\geq 1} \lambda(n) n^{iv}V\left(\frac{n}{N}\right)\\ & e\left(\frac{n\bar{a}}{q} - \frac{nx}{aq} \right) N^{1-i(t+v)} \sum_{m\equiv\bar{a}\mod q} U^\dagger\left(\frac{N(ma-x)}{aq}, 1-i(t+v)\right) dv dx.
\end{split}
\end{equation}

We can have $m=0$ only when $q=1$, in which case, $N(ma-x)/aq \ll N/Qq$, so its contribution to the sum will be negligible (as soon as $Q$ has size). 

For $m\neq0$, we have $N(ma-x)/aq \asymp N|m|/q$. Bounds on $U^\dagger$ give
\begin{equation}
U^\dagger \left(\frac{N(ma-x)}{aq}, 1-i(t+v)\right) \ll_j \left(\frac{1+|t+v|}{N|m|q^{-1}}\right)^j
\end{equation}
Thus we get arbitrary saving for $|m|>(1+|t+v|)q/N$. If we make sure $v<t$, that is $K<t$,  we'll have arbitrary saving for $|m|\gg qt^{1+\epsilon}/N$. Noting the condition $m\equiv \bar{a} \mod q$ and rearranging the sums in $S^+(N)$, 

\begin{equation}\label{poisson}
\begin{split}
S^+(N)= & \frac{N}{K}\int_0^1\int_{\R} N^{-i(t+v)}V\left(\frac{v}{K}\right)\sum_{1\leq q\leq Q}\underset{(m,q)=1}{\sum_{1\leq |m|\ll \frac{qt^{1+\epsilon}}{N}}}\frac{1}{aq} U^\dagger\left(\frac{N(ma-x)}{aq},1-i(t+v)\right) \\ & \sum_{n\geq1} \lambda(n) n^{iv}V\left(\frac{n}{N}\right) e\left(\frac{nm}{q}-\frac{nx}{aq}\right) dv dx 
\end{split}
\end{equation}
where $a\in (Q,q+Q]$ is the unique multiplicative inverse of $m\mod q$.
\begin{rem} Trivial bound here gives $S^+(N) \ll Nt^{1+\epsilon}$. We need to save $t$ and a bit more.
\end{rem}

We next split the $q-$sum into dyadic segments $(C,2C]$ 

\begin{equation*}
S^+(N) = \frac{N}{K} \sum_{1\leq C\leq Q} S(N,C)
\end{equation*}
where

\begin{equation}
\begin{split}
S(N,C) = & \int_0^1\int_{\R} N^{-i(t+v)}V\left(\frac{v}{K}\right)\sum_{C< q\leq 2C}\underset{(m,q)=1}{\sum_{1\leq |m|\ll \frac{qt^{1+\epsilon}}{N}}}\frac{1}{aq} U^\dagger\left(\frac{N(ma-x)}{aq},1-i(t+v)\right) \\ & \sum_{n\geq1} \lambda(n) n^{iv}V\left(\frac{n}{N}\right) e\left(\frac{nm}{q}-\frac{nx}{aq}\right) dv dx.
\end{split}
\end{equation}

\subsection{Voronoi summation to the $n$-sum} 
Applying Lemma \ref{voronoi} to the $n$-sum gets us
\begin{equation}
\begin{split}
\sum_{n\geq1}\lambda_f(n)e\left(n\frac{m}{q}\right)F(n) = & \frac{\pi i^k}{q} \sum_{n\geq1}\lambda_{f}(n)e\left(-n\frac{a}{q}\right)\int_0^{\infty} y^{iv} V\left(\frac{y}{N}\right) e\left(\frac{-xy}{aq}\right) \\& \times \frac{1}{2\pi i} \int_{(\sigma)} \left(\frac{2\pi\sqrt{ny}}{q}\right)^{-s} \frac{\Gamma(s/2 + (k-1)/2)}{\Gamma(1-s/2 + (k-1)/2)} ds dy 
\end{split}
\end{equation}
where $F(y) = y^{iv}V(\frac{y}{N})e(\frac{-xy}{aq})$. We want to be able to interchange integrals. For this, we use the complex Stirling approximation
\begin{equation*}
|\Gamma(z)| = \sqrt{2\pi}e^{-\sigma}|z|^{\sigma-1/2}e^{-\tau \arg(z)}\left(1+O\left(\frac{1}{|z|}\right)\right)
\end{equation*}
for $\arg(z)<\pi$ and $|z|\rightarrow\infty$. For
\begin{equation*}
\gamma(s) = (2\pi)^{-s}\frac{\Gamma(s/2 + (k-1)/2)}{\Gamma(1-s/2+(k-1)/2)}
\end{equation*}
we have
\begin{equation*}
|\gamma(s)| \sim (2\pi)^{-\sigma}e^{1-\sigma}|\tau|^{\sigma-1} \text{\quad as $|\tau| \rightarrow \infty$}
\end{equation*}
Looking at the pole free regions of the $\Gamma-$factors in the definition of $\gamma(s)$, we get 
\begin{equation}
|\gamma(s)| \ll 1 + |\tau|^{\sigma-1} \text{\quad for $\sigma>1-k$}
\end{equation}
We cannot apply Fubini theorem to interchange integrals right away since the integral is not absolutely convergent for $0<\sigma<1$. But if we assume that $k>1$, we can shift the integral to the line $\sigma= -1/2$ without picking any residues and the integral would be absolutely convergent, allowing us to apply Fubini and interchange integrals. 

\begin{equation*}
\begin{split}
\sum_{n\geq1}\lambda_f(n)e\left(n\frac{m}{q}\right)F(n) = & \frac{\pi i^k}{q} \sum_{n\geq1}\lambda_{f}(n)e\left(-n\frac{a}{q}\right)\frac{1}{2\pi i} \int_{(-1/2)} \left(\frac{\sqrt{n}}{q}\right)^{-s} \gamma(s) \\& \times \int_0^{\infty} y^{-s/2+iv} V\left(\frac{y}{N}\right) e\left(\frac{-xy}{aq}\right)  dy ds \\= & \frac{\pi i^kN^{1+iv}}{q} \sum_{n\geq1}\lambda_{f}(n)e\left(-n\frac{a}{q}\right)\frac{1}{2\pi i} \int_{(-1/2)} \left(\frac{\sqrt{nN}}{q}\right)^{-s} \gamma(s) \\& \times \int_0^{\infty} y^{-s/2+iv} V(y) e\left(\frac{-xN}{aq}y\right)  dy ds \\= & \frac{\pi i^kN^{1+iv}}{q} \sum_{n\geq1}\lambda_{f}(n)e\left(-n\frac{a}{q}\right)\frac{1}{2\pi i} \int_{(-1/2)} \left(\frac{\sqrt{nN}}{q}\right)^{-s} \gamma(s) \\& \times V^\dagger\left(\frac{xN}{aq}, 1-s/2+iv \right) ds
\end{split}
\end{equation*}

The bound on $V^\dagger$ gives
\begin{equation}
V^\dagger\left(\frac{xN}{aq}, 1-s/2+iv \right) \ll_j \min \left\lbrace 1 , \left(\frac{1 + |Nx/aq|}{|v-\tau/2|}\right)^j \right\rbrace
\end{equation}
We can therefore shift the integral from $\sigma = -1/2$ to $\sigma=M$ for any large $M$ by choosing $j=M+1$ (which kills the growth of $\gamma(s)$). We'll thus get saving for large $n$.

\begin{rem}\label{conductorlowering}
Using the above bound on $V^{\dagger}$, we get
\begin{equation*}
\left(\frac{\sqrt{nN}}{q}\right)^{-s}\gamma(s)V^{\dagger}\left(\frac{Nx}{aq},1-\frac{s}{2}+iv \right) \ll_j \left(\frac{\sqrt{nN}}{q}\right)^{-M} (1+|\tau|^{M-1})\min\left\lbrace 1, \left(\frac{1 + |Nx/aq|}{|v-\tau/2|}\right)^j \right\rbrace
\end{equation*}
Since $v\asymp K$, the better bound on $V^{\dagger}$ would be $O(1)$ when $|\tau|\leq 8K$. In that case,
\begin{equation*}
\begin{split}
\int_{|\tau|\leq 8K} \left(\frac{\sqrt{nN}}{q}\right)^{-s}\gamma(s)V^{\dagger}\left(\frac{Nx}{aq},1-\frac{s}{2}+iv \right) & \ll \int_{|\tau|\leq 8K} \left(\frac{\sqrt{nN}}{q}\right)^{-M} |\tau|^{M-1} d\tau \\ & \ll \left(\frac{\sqrt{nN}}{qK}\right)^{-M}
\end{split}
\end{equation*}
We'll thus get arbitrary saving for $n\gg Q^2K^2t^{\epsilon}/N$. On the other hand, when $|\tau|>8K$, we have the bound $V^{\dagger}\ll_j(N/aq|\tau|)^j$. Taking $j=M+1$,
\begin{equation*}
\begin{split}
\int_{|\tau|>8K} \left(\frac{\sqrt{nN}}{q}\right)^{-s}\gamma(s)V^{\dagger}\left(\frac{Nx}{aq},1-\frac{s}{2}+iv \right) & \ll \int_{|\tau|>8K} \left(\frac{\sqrt{nN}}{q}\right)^{-M} |\tau|^{M-1} \left(\frac{N}{aq|\tau|}\right)^{M+1} d\tau \\ & \ll \left(\frac{\sqrt{nN}}{q}\right)^{-M} \left(\frac{N}{aq}\right)^{M+1} \\ & = \left(\frac{an^{1/2}}{N^{1/2}}\right)^{-M} \left(\frac{N}{aq}\right)^2
\end{split}
\end{equation*}
We'll thus get arbitrary saving for $n\gg Nt^{\epsilon}/Q^2$. It makes sense to choose $Q$ so that the two bounds on $n$ are equal. Therefore set $Q=(N/K)^{1/2}$. We'll get arbitrary saving for $n\gg Kt^{\epsilon}$.
\end{rem}

For smaller values of $n$, we take $\sigma=1$. Note that the $\gamma$ factor will then be bounded.

\begin{equation}\label{sigma1}
\begin{split}
\sum_{n\geq1}\lambda_f(n)e\left(n\frac{m}{q}\right)F(n) = & \pi i^kN^{1/2+iv} \sum_{n\ll Q^2K^2/N}\frac{\lambda_{f}(n)}{n^{1/2}}e\left(-n\frac{a}{q}\right)\frac{1}{2\pi} \int_{\R} \left(\frac{\sqrt{nN}}{q}\right)^{-i\tau} \\& \times \gamma(1+i\tau) V^\dagger\left(\frac{xN}{aq}, 1/2-i\tau/2+iv \right) d\tau
\end{split}
\end{equation}

Assuming $K\ll t^{1-\epsilon}$, we get arbitrary saving for $|\tau| > Nt^\epsilon/QC$ due to bounds on $V^\dagger$. Thus we can restrict the integral to $\tau \in [-Nt^\epsilon/QC, Nt^\epsilon/QC]$ by defining a smooth partition of unity on this set. Let $W_J$  for $J\in\mathcal{J}$ be smooth bump functions satisfying $x^lW_J^{(l)}\ll_l 1$ for all $l\geq0$. For $J=0$, let the support of $W_0$ be in $[-1,1]$ and for $J>0$ (resp. $J<0$), let the support of $W_J$ be in $[J,4J/3]$ (resp $[4J/3,J]$). Finally, we require that
\begin{equation*}
\sum_{J\in \mathcal{J}} W_J(x) = 1 \text{\quad for \quad $x \in [-Nt^\epsilon/QC, Nt^\epsilon/QC]$}
\end{equation*}
The precise definition of the functions $W_J$ will not be needed. We note that we need only $O(\log(t))$ such $J\in\mathcal{J}$. We can write the integral appearing in Voronoi summation as

\begin{equation*}
\begin{split}
& \int_{\R} \left(\frac{\sqrt{nN}}{q}\right)^{-i\tau} \gamma(1+i\tau) V^\dagger\left(\frac{xN}{aq}, 1/2-i\tau/2+iv \right) d\tau =\\ & \sum_{J\in \mathcal{J}} \int_{\R} \left(\frac{\sqrt{nN}}{q}\right)^{-i\tau} \gamma(1+i\tau) V^\dagger\left(\frac{xN}{aq}, 1/2-i\tau/2+iv \right) W_J(\tau) d\tau + O(t^{-20150})
\end{split}
\end{equation*}

Combining everything, we write $S(N,C)$ as
\begin{equation}
\begin{split}
S(N,C) = \frac{i^kN^{1/2-it}K}{2} & \sum_{J\in\mathcal{J}}\sum_{n\ll Q^2K^2/N} \frac{\lambda_{f}(n)}{n^{1/2}} \sum_{C<q\leq2C} \underset{1\leq|m|\ll \frac{qt^{1+\epsilon}}{N}}{\sum_{(m,q)=1}} e\left(\frac{-na}{q}\right) \frac{1}{aq} \int_{\R} \left(\frac{\sqrt{nN}}{q}\right)^{-i\tau} \\& \times \gamma\left(1 + i\tau\right) W_J(\tau) \mathcal{I^{**}}(q,m,\tau) d\tau + O(t^{-2015})
\end{split}
\end{equation}
where

\begin{equation*}
\mathcal{I^{**}}(q,m,\tau)= \int_0^1 \int_{\R} V(v) U^{\dagger}\left( \frac{N(ma-x)}{aq}, 1-i(t+Kv) \right) V^{\dagger}\left( \frac{Nx}{aq}, \frac{1}{2} - \frac{i\tau}{2} + iKv \right) dv dx
\end{equation*}

\begin{rem}
We can trivially bound $I^{**}(q,m,\tau)$ by $O(1)$, and the $\tau$-integral is over the interval $[-Nt^\varepsilon/QC, Nt^\varepsilon/QC]$. Trivial bound on $S(N,C)$ will imply $S(N,C)\ll K^{5/2}t^{1+\varepsilon}/N^{1/2}$. So $S(N)\ll N^{1/2}K^{3/2}t^{1+\varepsilon}$. We need to save $N^{1/2}K^{3/2}$ and a bit more.
\end{rem}

\section{Analysis of the integrals}

We next analyze the integral $\mathcal{I}^{**}(q,m,\tau)$. Application of Lemma \ref{dagger} to $U^{\dagger}$ gives us
\begin{equation*}
U^{\dagger} = \frac{e^{i\pi/4}(t+Kv)^{1/2}aq}{(2\pi)^{1/2}N(x-ma)}U\left( \frac{(t+Kv)aq}{2\pi N(x-ma)} \right) \left( \frac{(t+Kv)aq}{2\pi eN(x-ma)} \right)^{-i(t+Kv)} + O(t^{-3/2}).
\end{equation*}
Therefore,
\begin{equation*}
\begin{split}
\mathcal{I^{**}}(\tau)= \frac{c_1aq}{N} \int_0^1 \int_{\R} & V(v) V^{\dagger}\left( \frac{Nx}{aq}, \frac{1}{2} - \frac{i\tau}{2} + iKv \right) \frac{(t+Kv)^{1/2}}{(x-ma)}U\left( \frac{(t+Kv)aq}{2\pi N(x-ma)} \right) \\& \times \left( \frac{(t+Kv)aq}{2\pi eN(x-ma)} \right)^{-i(t+Kv)} dv dx + O(t^{-3/2+\epsilon})
\end{split}
\end{equation*}
where $c_1 = e^{i\pi/4}/\sqrt{2\pi}$. We next apply Lemma \ref{dagger} to $V^{\dagger}$.
\begin{equation*}
\begin{split}
V^{\dagger} = \frac{2\sqrt{\pi} e^{-i\pi/4}}{(4\pi)^{1/2}} \left(\frac{aq}{Nx}\right)^{1/2} & V\left(\frac{(2Kv-\tau)aq}{4\pi Nx}\right) \left(\frac{(2Kv-\tau)aq}{4\pi eNx}\right)^{i(Kv-\tau/2)} \\ & + O\left(\min\left\lbrace \left(\frac{aq}{Nx}\right)^{3/2}, \frac{1}{|\tau/2-Kv|^{3/2}} \right\rbrace \right).
\end{split}
\end{equation*}

The integral then becomes
\begin{equation}\label{I_tau}
\begin{split}
\mathcal{I^{**}}(q,m,\tau)= c_2\left(\frac{aq}{N}\right)^{3/2} \int_0^1 \int_{\R} & V(v) \left(\frac{1}{x}\right)^{1/2} V\left(\frac{(2Kv-\tau)aq}{4\pi Nx}\right) \left(\frac{(2Kv-\tau)aq}{4\pi eNx}\right)^{i(Kv-\tau/2)} \\& \times \frac{(t+Kv)^{1/2}}{(x-ma)}U\left( \frac{(t+Kv)aq}{2\pi N(x-ma)} \right) \left( \frac{(t+Kv)aq}{2\pi eN(x-ma)} \right)^{-i(t+Kv)} dv dx \\& + O(E^{**} + t^{-3/2+\epsilon})
\end{split}
\end{equation}
with $c_2= 1/(2\pi)^{1/2}$ and since $uU(u)\ll 1$, 
\begin{equation*}
E^{**} = \frac{1}{t^{1/2}}\int_0^1\int_1^2 \min\left\lbrace \left(\frac{aq}{Nx}\right)^{3/2}, \frac{1}{|\tau/2-Kv|^{3/2}} \right\rbrace dv dx
\end{equation*}
(We note that more generally $u^jU(u)\ll_j 1$, but using this does not improve the error term.)
\subsection{Analysis of the error term $E^{**}$}
The first term is smaller than the second if and only if
\begin{equation*}
\frac{\tau}{2K} - \frac{Nx}{aqK} < v < \frac{\tau}{2K} + \frac{Nx}{aqK}.
\end{equation*}

If $|\tau|\geq 10K$, this interval does not intersect $[1,2]$ unless $Nx/aq \asymp |\tau|$. For this, we use the trivial bound $O(1)$ for the inner integral over $v$. And if $|\tau|<10K$, the inner integral is bounded by the length of the interval, which is $2Nx/aqK$. Hence the contribution where the first term is smaller than the second is of the order
\begin{equation*}
\frac{1}{t^{1/2}}\int_0^1 \left(\frac{aq}{Nx}\right)^{1/2} \frac{1}{K}\textbf{1}_{|\tau|<10K} dx + \frac{1}{t^{1/2}}\int_0^1 \left(\frac{aq}{Nx}\right)^{1/2} \frac{1}{|\tau|}\textbf{1}_{|\tau|\geq10K} dx.
\end{equation*}
This is bounded by
\begin{equation*}
O\left(\frac{Q}{t^{1/2}N^{1/2}K}\min\left\lbrace 1, \frac{10K}{|\tau|} \right\rbrace t^\epsilon \right).
\end{equation*}

Next we estimate the contribution to $E^{**}$ when the second term is smaller.  This would be
\begin{equation*}
\begin{split}
\frac{1}{t^{1/2}}\underset{|\tau/2 - Kv| > Nx/aq}{\int_0^1\int_1^2} \frac{1}{|\tau - Kv|^{3/2}} dv dv & \ll \frac{1}{t^{1/2}} \int_0^1 \left(\frac{aq}{Nx}\right)^{1/2+\epsilon} \int_1^2 \frac{1}{|\tau/2 - Kv|^{1-\epsilon}} dv dx \\ & \ll t^{\epsilon} \frac{Q}{t^{1/2}N^{1/2}K}\min\left\lbrace 1, \frac{10K}{|\tau|} \right\rbrace.
\end{split}
\end{equation*}

The total error term therefore is 
\begin{equation}\label{E1}
E^{**} + t^{-3/2 + \epsilon} \ll t^{\epsilon} \frac{Q}{t^{1/2}N^{1/2}K}\min\left\lbrace 1, \frac{10K}{|\tau|} \right\rbrace + t^{-3/2+\epsilon},
\end{equation}
and we can write
\begin{equation}
\begin{split}
\mathcal{I^{**}}(q,m,\tau)= c_2\left(\frac{aq}{N}\right)^{3/2} \int_0^1 \int_{\R} & V(v) \left(\frac{1}{x}\right)^{1/2} V\left(\frac{(2Kv-\tau)aq}{4\pi Nx}\right) \left(\frac{(2Kv-\tau)aq}{4\pi eNx}\right)^{i(Kv-\tau/2)} \\& \times \frac{(t+Kv)^{1/2}}{(x-ma)}U\left( \frac{(t+Kv)aq}{2\pi N(x-ma)} \right) \left( \frac{(t+Kv)aq}{2\pi eN(x-ma)} \right)^{-i(t+Kv)} dv dx \\& + O\left(\frac{t^{\epsilon}}{t^{1/2}K^{3/2}}\min\left\lbrace 1, \frac{10K}{|\tau|}\right\rbrace + t^{-3/2+\epsilon}\right).
\end{split}
\end{equation}

\begin{rem}
The error term in the above estimate for $I^{**}$ saves a further $t^{1/2}K^{3/2}$. The main term saves $K^{1/2}t^{1/2}$. So we need to save $K$ and a bit more. Note that at this point $K$ seems to be hurting us rather than helping us. Moreover, if $K$ had no size, we would get the bound $S(N)\ll N^{1+\varepsilon}$, which would get us the convexity bound.
\end{rem}

\subsection{Analysis of integral over $v$} The integral is given by 
\begin{equation*}
\begin{split}
I_1 = c_2\left(\frac{aq}{N}\right)^{3/2} \int_0^1 \int_{\R} & V(v) \left(\frac{1}{x}\right)^{1/2} V\left(\frac{(2Kv-\tau)aq}{4\pi Nx}\right) \left(\frac{(2Kv-\tau)aq}{4\pi eNx}\right)^{i(Kv-\tau/2)} \\& \times \frac{(t+Kv)^{1/2}}{(x-ma)}U\left( \frac{(t+Kv)aq}{2\pi N(x-ma)} \right) \left( \frac{(t+Kv)aq}{2\pi eN(x-ma)} \right)^{-i(t+Kv)} dv dx
\end{split}
\end{equation*}
\\
Due to the argument of $U$, the integral vanishes if $m>0$. Trivial estimate gives
\begin{equation*}
\begin{split}
I_1 \ll \left(\frac{aq}{N}\right)^{3/2} \int_0^1 \int_{\R} & \frac{(t+Kv)^{1/2}}{x^{1/2}(x-ma)} V(v) V\left(\frac{(2Kv-\tau)aq}{4\pi Nx}\right) U\left( \frac{(t+Kv)aq}{2\pi N(x-ma)} \right) dv dx
\end{split}
\end{equation*}
\\
The length of the integral over $v$ is restricted due to the weight functions, respectively given by $1, -Nm/Kq$ and $Nx/aqK$. $Nx/aqK < -Nm/Kq$, so we can restrict the length of integral over $v$ to $Nx/aqK$. We restrict the integral over $x$ to $[0,1/K]$ and estimate the resulting integral trivially.

\begin{equation*}
\begin{split}
& \left(\frac{aq}{N}\right)^{3/2} \int_0^{1/K} \int_{\R} \frac{(t+Kv)^{1/2}}{x^{1/2}(x-ma)} V(v) V\left(\frac{(2Kv-\tau)aq}{4\pi Nx}\right) U\left( \frac{(t+Kv)aq}{2\pi N(x-ma)} \right) dv dx \\ & \ll \left(\frac{aq}{N}\right)^{1/2} \frac{1}{t^{1/2}} \int_0^{1/K} \frac{1}{x^{1/2}}\frac{Nx}{aqK} dx \\ &\ll \frac{1}{t^{1/2}K^{3/2 + 1}} \left(\frac{N}{aq}\right)^{1/2} = E
\end{split}
\end{equation*}
\\\\
We write $I_1(\tau) = I_2(\tau) + O(E)$, where $I_2(\tau)$ is
\begin{equation*}
\begin{split}
I_2 = c_2\frac{1}{t^{1/2}}\left(\frac{aq}{N}\right)^{3/2} \int_{1/K}^1 \int_{\R} & \frac{t^{1/2}(t+Kv)^{1/2}}{(x-ma)x^{1/2}} V(v) V\left(\frac{(2Kv-\tau)aq}{4\pi Nx}\right) \left(\frac{(2Kv-\tau)aq}{4\pi eNx}\right)^{i(Kv-\tau/2)} \\& \times U\left( \frac{(t+Kv)aq}{2\pi N(x-ma)} \right) \left( \frac{(t+Kv)aq}{2\pi eN(x-ma)} \right)^{-i(t+Kv)} dv dx
\end{split}
\end{equation*}
where an extra $t^{1/2}$ is multiplied to balance the size of the function. Set
\begin{equation*}
f(v) = -\frac{t+Kv}{2\pi}\log \left( \frac{(t+Kv)aq}{2\pi eN(x-ma)} \right) + \frac{2Kv - \tau}{4\pi} \log \left( \frac{(2Kv - \tau)aq}{4\pi eNx} \right)
\end{equation*}
and
\begin{equation*}
g(v) = \frac{t^{1/2}(t+Kv)^{1/2}aq}{N(x-ma)} V(v) V\left(\frac{(2Kv-\tau)aq}{4\pi Nx}\right) U\left( \frac{(t+Kv)aq}{2\pi N(x-ma)} \right)
\end{equation*}
So that
\begin{equation*}
I_2 = c_2\frac{1}{t^{1/2}}\left(\frac{aq}{N}\right)^{1/2}\int_{1/K}^1 \frac{1}{x^{1/2}} \int_{\R} g(v)e(f(v))dv dx
\end{equation*}
Then
\begin{equation*}
f'(v) = -\frac{K}{2\pi} \log \left( \frac{2(t+Kv)x}{(2Kv-\tau)(x-ma)} \right), \quad f^{(j)}(v) = -\frac{(j-1)!(-K)^j}{2\pi(t+Kv)^{j-1}} + \frac{(j-1)!(-2K)^j}{4\pi(2Kv-\tau)^{j-1}}
\end{equation*}
The stationary phase is given by
\begin{equation*}
v_0 = -\frac{(2t+\tau)x - \tau ma}{2Kma}
\end{equation*}

In support of the integral, we have
\begin{equation*}
f^{(j)} (v) \asymp \frac{Nx}{aq} \left( \frac{Kaq}{Nx} \right)^j
\end{equation*}
for $j\geq 2$, and for $j\geq 0$
\begin{equation*}
g^{(j)}(v) \ll \left( 1 + \frac{Kaq}{Nx} \right)^j
\end{equation*}

We shall apply the sharp version of stationary phase method due to Huxley\cite{huxley} (as given in Lemma 3 of Munshi\cite{munshi}):

We can write 
\begin{equation*}
f'(v) = \frac{K}{2\pi}\log\left(1+\frac{K(v_0-v)}{(t+Kv)}\right) - \frac{K}{2\pi}\log\left(1 + \frac{2K(v_0-v)}{(2Kv-\tau)}\right)
\end{equation*}

In the support of the integral, we have $0\leq 2Kv-\tau \ll N/aq \ll t^{1+\epsilon}/Q$ (since $N/t^{1+\epsilon}<q$ and $a \asymp Q$). Therefore
\begin{equation*}
f''(v) = -\frac{K^2}{2\pi(t+Kv)} + \frac{K^2}{2\pi(Kv-\tau/2)}
\end{equation*}
is positive on the support of the integral for large enough $t$. So $f'$ changes sign from negative to positive at $v_0$. Support of the integral is contained in $[1,2]$ due the weight function $V(v)$. If $v_0\notin [0.5,2.5]$, then $v_0$ is not in the support of the integral and $|v_0-v|>0.5$.
In the support of the integral, we will have
\begin{equation*}
|f'(v)|\gg K^{1-\epsilon}\min\left\lbrace1,\frac{Kaq}{Nx}\right\rbrace
\end{equation*}
Applying the first statement of Lemma (\ref{stationary}) with
\begin{equation*}
\Theta_f = \frac{Nx}{aq}, \quad \Omega_f = \frac{Nx}{Kaq}, \quad \Omega_g = \min\left\lbrace1,\frac{Nx}{Kaq}\right\rbrace, \quad \Lambda = K^{1-\epsilon}\min\left\lbrace1,\frac{Kaq}{Nx}\right\rbrace
\end{equation*}
we obtain the bound
\begin{equation}\label{b1}
\int_{\R} g(x) e(f(x)) dx \ll \frac{\Theta_f}{\Omega_f^2\Lambda^3}\left(1 + \frac{\Omega_f}{\Omega_g} + \frac{\Omega_f^2}{\Omega_g^2}\frac{\Lambda}{\Theta_f/\Omega_f}\right)t^{\epsilon}
\end{equation}
On the other hand, if $v_0 \in [0.5,2.5]$, then treating the integral as one over the finite range $[0.1,4]$ (so that $\kappa>0.4$) and applying the second part of Lemma (\ref{stationary}), we get
\begin{equation}\label{b2}
I = \frac{g(x_0)e(f(x_0) + 1/8)}{\sqrt{f''(x_0)}} + O\left(\left(\frac{\Omega_f^4}{\Theta_f^2} + \frac{\Omega_f}{\Theta_f^{3/2}} + \frac{\Omega_f^3}{\Theta_f^{3/2}\Omega_g^2}\right)t^{\epsilon}\right)
\end{equation} 

For the range $x\in[1/K,1]$, we use the bound in lemma (\ref{stationary}). In the case there is no stationary phase, we will use the first statement of lemma (\ref{stationary}). We have,
\begin{equation}
\Theta_f=\frac{Nx}{aq},\quad\Omega_f=\frac{Nx}{aqK},\quad\Lambda=K^{1-\varepsilon}\min\left\lbrace1,\frac{Kaq}{Nx}\right\rbrace,\quad\Omega_g=\min\left\lbrace1,\frac{Nx}{aqK}\right\rbrace.
\end{equation}

Next is the contribution of $x\in[1/K,1]$ when there is no stationary phase. When $x<aqK/N$, $\Lambda=K$ and $\Omega_g=\Omega_f$. In that case, the contribution is 
\begin{equation*}
\left(\frac{2\pi aq}{Nt}\right)^{1/2}\int_{1/K}^{\max\{\frac{1}{K},\frac{Kaq}{N}\}}\frac{1}{x^{1/2}}\frac{aq}{NKx}dx\ll \frac{1}{t^{1/2}K^{2}}.
\end{equation*}
This is always smaller than the contribution of the bound $E$. When $x>aqK/N$, $\Lambda=K^2aq/Nx$ and $\Omega_g=1$. In that case, the contribution is $1/K^3t^{1/2}$, which is better than above. We next calculate the contribution of the error term when there is a stationary phase. For that we have $\kappa>0.4$. One can calculate that for both $x<aqK/N$ and $x>aqK/N$, the contribution is $1/K^2t^{1/2}$. 

With all of this, we summarize the analysis in the following Lemma. Let
\begin{equation}\label{BCTau}
B(C,\tau)= \frac{t^\varepsilon}{t^{1/2}K^{3/2}}\min\left\lbrace1,\frac{10K}{|\tau|}\right\rbrace + \frac{1}{t^{1/2}K^{5/2}}\left(\frac{N}{QC}\right)^{1/2}.
\end{equation}

Note that,
\begin{equation}\label{intB}
\int_{-Nt^{\varepsilon}/QC}^{Nt^{\varepsilon}/QC}B(C,\tau)d\tau\ll \frac{K}{t^{1/2}K^{3/2}}+\frac{1}{t^{1/2}K^{5/2}}\left(\frac{N}{QC}\right)^{3/2}.
\end{equation}

Putting everything together, we have
\begin{lem}\label{I**}
Suppose $C<q\leq 2C$, with $1\ll C\leq (N/K)^{1/2}$ and $K$ satisfies $1\leq K\ll t^{1-\epsilon}$. Suppose $t>2$ and $|\tau|\ll N^{1/2}K^{1/2}t^{\epsilon}$. We have
\begin{equation*}
\mathcal{I}^{**}(q,m,\tau) = \mathcal{I}_1(q,m,\tau) + \mathcal{I}_2(q,m,\tau)
\end{equation*}
where
\begin{equation*}
\mathcal{I}_1(q,m,\tau) = \frac{c_4}{(t+\tau/2)^{1/2}K}\left(-\frac{(t+\tau/2)q}{2\pi eNm}\right)^{3/2-i(t+\tau/2)}V\left(-\frac{(t+\tau/2)q}{2\pi Nm}\right)\int_0^1 V\left(\frac{\tau}{2K} - \frac{(t+\tau/2)x}{Kma}\right) dx
\end{equation*}
for some absolute constant $c_4$ and
\begin{equation*}
\mathcal{I}_2(q,m,\tau) := \mathcal{I}^{**}(q,m,\tau) -  \mathcal{I}_1(q,m,\tau) = O(B(C,\tau)t^{\epsilon})
\end{equation*}
with $B(C,\tau)$ as defined in (\ref{BCTau}).
\end{lem}
Consequently, we have the following decomposition of $S(N,C)$.

\begin{lem}
\begin{equation*}
S(N,C) = \sum_{J\in\mathcal{J}}\{S_{1, J}(N,C) + S_{2,J}(N, C)\} + O(t^{-2015})
\end{equation*}
where
\begin{equation*}
S_{l, J}(N,C) = \frac{i^kN^{1/2-it}K}{2} \sum_{n\ll Q^2K^2/N} \frac{\lambda_{f}(n)}{n^{1/2}} \sum_{C<q\leq2C} \underset{1\leq|m|\ll \frac{qt^{1+\epsilon}}{N}}{\sum_{(m,q)=1}}e\left(\frac{-na}{q}\right) \frac{1}{aq}\mathcal{I}_{l,J}(q,m,n)
\end{equation*}
and
\begin{equation*}
\mathcal{I}_{l,J}(q,m,n) = \int_{\R} \left(\frac{\sqrt{nN}}{q}\right)^{-i\tau} \gamma\left(1 + i\tau\right) W_J(\tau) \mathcal{I}_l(q,m,\tau) d\tau
\end{equation*}
with $\mathcal{I}_l(q,m,\tau)$ as defined in the previous lemma.
\end{lem}

\begin{rem}
The saving due to $\mathcal{I}_1(q,m,\tau)$ is still $t^{1/2}K^{1/2}$, same as the main term before this analysis. The saving due to $\mathcal{I}_2(q,m,\tau)$ is $t^{1/2}K^{9/4}/N^{1/4}$. In all, we need to save $\max\{K,t^{1/4}/K^{3/4}\}$ and a bit more.
\end{rem}

\section{Application of Cauchy and Poisson summation- I}

In this section, we will estimate
\begin{equation*}
S_2(N,C) := \sum_{J\in\mathcal{J}} S_{2,J}(N,C)
\end{equation*}
Here, we'll not apply any cancellation over the $\tau$-integral. Dividing the $n$-sum into dyadic segments and using the bound $\gamma(1+i\tau)\ll 1$,  we get
\begin{equation}
S_2(N,C) \ll t^{\epsilon}N^{1/2}K\int_{-\frac{(NK)^{1/2}t^{\epsilon}}{C}}^{\frac{(NK)^{1/2}t^{\epsilon}}{C}}\underset{dyadic}{\sum_{1\leq L \ll Kt^{\epsilon}}}\sum_n \frac{|\lambda_f(n)|}{n^{1/2}}U\left(\frac{n}{L}\right)\bigg|\underset{1\leq|m|\ll \frac{qt^{1+\epsilon}}{N}}{\sum_{C<q\leq 2C}\sum_{(m,q)=1}}e\left(\frac{-na}{q}\right)\frac{1}{aq^{1-i\tau}}\mathcal{I}_2(q,m,\tau)\bigg|d\tau
\end{equation}
Applying Cauchy to the $n-$sum and using the Ramanujan bound on average (Lemma \ref{ramonavg}), we get
\begin{equation}
S_2(N,C) \ll t^{\epsilon}N^{1/2}K\int_{-\frac{(NK)^{1/2}t^{\epsilon}}{C}}^{\frac{(NK)^{1/2}t^{\epsilon}}{C}}\underset{dyadic}{\sum_{1\leq L \ll Kt^{\epsilon}}} L^{1/2}[S_2(N,C,L,\tau)]^{1/2}d\tau
\end{equation}
where
\begin{equation*}
\begin{split}
S_2(N,C,L,\tau) = & \sum_n\frac{1}{n}U\left(\frac{n}{L}\right) \mid\underset{1\leq|m|\ll \frac{qt^{1+\epsilon}}{N}}{\sum_{C<q\leq 2C}\sum_{(m,q)=1}}e\left(\frac{-na}{q}\right)\frac{1}{aq^{1-i\tau}}\mathcal{I}_2(q,m,\tau)\mid^2 \\ = & \sum_n \frac{1}{n}U\left(\frac{n}{L}\right) \underset{1\leq|m|\ll \frac{qt^{1+\epsilon}}{N}}{\sum_{C<q\leq 2C}\sum_{(m,q)=1}}e\left(\frac{-na}{q}\right)\frac{1}{aq^{1-i\tau}}\mathcal{I}_2(q,m,\tau) \\ & \times \underset{1\leq|m'|\ll \frac{q't^{1+\epsilon}}{N}}{\sum_{C<q'\leq 2C}\sum_{(m',q')=1}}e\left(\frac{na'}{q'}\right)\frac{1}{a'q'^{1+i\tau}}\overline{\mathcal{I}_2(q',m',\tau)} \\ = & \underset{1\leq|m|\ll \frac{qt^{1+\epsilon}}{N}}{\sum_{C<q\leq 2C}\sum_{(m,q)=1}} \underset{1\leq|m'|\ll \frac{q't^{1+\epsilon}}{N}}{\sum_{C<q'\leq 2C}\sum_{(m',q')=1}}\frac{1}{aq^{1-i\tau}} \frac{1}{a'q'^{1+i\tau}} \mathcal{I}_2(q,m,\tau) \overline{\mathcal{I}_2(q',m',\tau)} T
\end{split} 
\end{equation*}
where we set
\begin{equation*}
T = \sum_n \frac{1}{n} U\left(\frac{n}{L}\right)e\left(\frac{-na}{q}\right)e\left(\frac{na'}{q'}\right)
\end{equation*}
We break the $n-$sum modulo $qq'$ to get
\begin{equation*}
T = \sum_{\beta \mod qq'}e\left(\frac{\beta(a'q-aq')}{qq'}\right)\sum_{l\in\Z} \frac{1}{\beta + lqq'}U\left(\frac{\beta + lqq'}{L}\right)
\end{equation*}
Applying Poisson summation formula to $l-$sum,
\begin{equation*}
T = \sum_{\beta \mod qq'}e\left(\frac{\beta(a'q-aq')}{qq'}\right)\sum_{n\in\Z} \int_{\R}\frac{1}{\beta + yqq'}U\left(\frac{\beta + yqq'}{L}\right)e(-ny) dy
\end{equation*}
Change variables $w = (\beta + yqq')/L$ to get
\begin{equation*}
T = \frac{1}{qq'} \sum_{\beta\mod qq'} e\left(\frac{\beta(a'q-aq')}{qq'}\right)\sum_{n\in\Z} e\left(\frac{n\beta}{qq'}\right) \int_{\R}\frac{1}{w}U(w)e\left(\frac{-nLw}{qq'}\right) dw
\end{equation*}
Integration by parts will give arbitrary saving for $n\gg C^2t^{\epsilon}/L$. Thus,
\begin{equation*}
T = \frac{1}{qq'} \sum_{n\ll \frac{C^2t^\epsilon}{L}} \left[\sum_{\beta\mod qq'}e\left(\frac{\beta(a'q-aq')}{qq'}\right)e\left(\frac{n\beta}{qq'}\right)\right]\int_{\R}\frac{1}{w}U(w) e\left(\frac{-nLw}{qq'}\right) dw + O(t^{-2015})
\end{equation*}
Plugging this in the expression for $S_2(N,C,L,\tau)$, we get

\begin{equation*}
S_2(N,C,L,\tau)\ll \frac{K}{NC^4}B(C,\tau)^2 \underset{1\leq|m|\ll \frac{qt^{1+\epsilon}}{N}}{\sum_{C<q\leq 2C}\sum_{(m,q)=1}} \underset{1\leq|m'|\ll \frac{q't^{1+\epsilon}}{N}}{\sum_{C<q'\leq 2C}\sum_{(m',q')=1}}\sum_{n\ll \frac{C^2t^\epsilon}{L}} |\mathfrak{C}| + O(t^{-2015})
\end{equation*}
where
\begin{equation*}
\mathfrak{C} = \sum_{\beta\mod qq'}e\left(\frac{\beta(a'q-aq')}{qq'}\right)e\left(\frac{n\beta}{qq'}\right)
\end{equation*}

Note that $\mathfrak{C} = qq' \delta (n\equiv aq' - a'q \quad \mod qq')$. Plugging that into the above expression and rearranging the sums, we get

\begin{lem}
\begin{equation*}
S_2(N,C,L,\tau)\ll \frac{K}{NC^2}B(C,\tau)^2 \sum_{n\ll \frac{C^2t^\epsilon}{L}} \underset{1\leq|m|\ll \frac{qt^{1+\epsilon}}{N}}{\sum_{C<q\leq 2C}\sum_{(m,q)=1}} \underset{1\leq|m'|\ll \frac{q't^{1+\epsilon}}{N}}{\sum_{C<q'\leq 2C}\sum_{(m',q')=1}} \delta(n\equiv aq' - a'q \quad \mod qq') + O(t^{-2015})
\end{equation*}
\end{lem}

We have to analyze the cases $n=0$ and $n\neq0$ separately. When $n=0$, the congruence condition above gives $q=q'$ and $a=a'$. For a given $m$, this fixes $m'$ up to a factor of $t^{1+\varepsilon}/N$. Moreover, in the case $Q^2<K$, that is, $K>N^{1/2}$, we'll have only $n=0$ for $L>C^2$. Therefore for $n\neq0$, we will let $L$ go up to $\min\{C^2,K\}$.

We note that the congruence condition implies $q|(n-aq')$ and $q'|(n+a'q)$. Since $a$ and $a'$ lie in an interval of length $q$, fixing $n, q$ and $q'$ fixes both $a$ and $a'$. That saves $q, q'$ in the $m, m'$-sums respectively.

\begin{rem}
We haven't used the conditions $(a,q)=1$ and $(a',q')=1$. But we can show that these conditions give us a saving of at most a power of $\log t$.
\end{rem}

Using $\mathcal{I}_2(q,m,\tau)\ll B(C,\tau)$, we get
\begin{equation*}
S_2(N,C,L,\tau) \ll t^{\epsilon} \frac{Kt^2 B(C,\tau)^2}{N^3}\bigg[\underbrace{1}_{n=0}+\underbrace{\frac{C^2}{L}}_{n\neq0}\bigg]
\end{equation*}
so that
\begin{equation*}
S_2(N,C,L,\tau)^{1/2} \ll t^{\epsilon} \frac{K^{1/2}tB(C,\tau)}{N^{3/2}}\bigg[1 + \frac{C}{L^{1/2}}\bigg]
\end{equation*}
Therefore,
\begin{equation*}
S_2(N,C)\ll t^{\epsilon}N^{1/2}K \int_{-\frac{(NK)^{1/2}t^{\epsilon}}{C}}^{\frac{(NK)^{1/2}t^{\epsilon}}{C}}\bigg[\underset{dyadic}{\sum_{1\leq L \ll Kt^{\epsilon}}} L^{1/2}.\frac{K^{1/2}tB(C,\tau)}{N^{3/2}}+\underset{dyadic}{\sum_{1\leq L \ll \min\{C^2,K\}t^{\epsilon}}}\frac{K^{1/2}tCB(C,\tau)}{N^{3/2}}\bigg]d\tau
\end{equation*}
If $K\geq N^{1/2}$, then the contribution of the second term is smaller than that of the first. So we neglect the second term. Summing over $L$, using (\ref{intB}) (and noting $N\asymp t^{1+\varepsilon}$), we get
\begin{equation*}
S_2(N,C)\ll t^{\epsilon}\frac{K^2t}{N}\left(\frac{1}{t^{1/2}K^{1/2}}+\frac{1}{t^{1/2}K^{5/2}}\left(\frac{N}{QC}\right)^{3/2}\right)
\end{equation*}
Multiplying by $N^{1/2}/K$ and summing over $C$ dyadically, 
\begin{equation}\label{S_2}
\begin{split}
\frac{S_2(N)}{N^{1/2}} & \ll t^{1/2+\varepsilon}\left(\frac{K^{1/2}}{N^{1/2}}+\frac{N^{1/4}}{K^{3/4}}\right)
\end{split}
\end{equation}
where $K\geq N^{1/2}$.

\section{Application of Cauchy and Poisson summation- II}
\begin{equation*}
\mathcal{I}_1(q,m,\tau) = \frac{c_4}{(t+\tau/2)^{1/2}K}\left(-\frac{(t+\tau/2)q}{2\pi eNm}\right)^{3/2-i(t+\tau/2)}V\left(-\frac{(t+\tau/2)q}{2\pi Nm}\right)\int_0^1 V\left(\frac{\tau}{2K} - \frac{(t+\tau/2)x}{Kma}\right) dx
\end{equation*}

\begin{equation*}
\begin{split}
S_{1,J}(N,C) = \frac{i^kN^{1/2-it}K}{2} & \underset{L dyadic}{\sum_{1\leq L\ll Kt^{\epsilon}}}\sum_n \frac{\lambda_f(n)}{n^{1/2}} U\left(\frac{n}{L}\right) \underset{\begin{subarray}{c} C<q\leq2C, (m,q)=1 \\ 1\leq |m|\ll qt^{1+\epsilon}/N\end{subarray}}{\sum\sum} e\left(\frac{-na}{q}\right) \frac{1}{aq} \\ & \times \int_{\R} \left(\frac{\sqrt{nN}}{q}\right)^{-i\tau}\gamma(1+i\tau)W_{J}(\tau)I_1(q,m,\tau) d\tau
\end{split}
\end{equation*}

Using the two, rearranging $q,m-$sums and integral, taking absolute values and using Cauchy, we get

\begin{equation}
|S_{1,J}(N,C)| \leq N^{1/2}K \underset{L-dyadic}{\sum_{1\leq L\ll Kt^{\epsilon}}} \left(\sum_n |\lambda_f(n)|^2 U\left(\frac{n}{L}\right)\right)^{1/2}[S_{1,J}(N,C,L)]^{1/2}
\end{equation}
where
\begin{equation}
S_{1,J}(N,C,L) = \sum_n \frac{1}{n} U\left(\frac{n}{L}\right) \mid \int_{\R} (\sqrt{nN})^{-i\tau}\gamma(1+i\tau) \underset{\begin{subarray}{c} C<q\leq2C, (m,q)=1 \\ 1\leq |m|\ll qt^{1+\epsilon}/N\end{subarray}}{\sum\sum}e\left(\frac{-na}{q}\right) \frac{1}{aq^{1-i\tau}}W_{J}(\tau)I_1(q,m,\tau) d\tau \mid^2
\end{equation}
Opening $|...|^2$ and rearranging sums and integrals
\begin{equation*}
\begin{split}
S_{1,J}(N,C,L) = \underset{\R^2}{\int\int} &(\sqrt{N})^{-i\tau+i\tau'}\gamma(1+i\tau)\overline{\gamma(1+i\tau')} W_J(\tau)W_J(\tau') \\ &\times\underset{\begin{subarray}{c} C<q\leq2C, (m,q)=1 \\ 1\leq |m|\ll qt^{1+\epsilon}/N\end{subarray}}{\sum\sum} \underset{\begin{subarray}{c} C<q'\leq2C, (m',q')=1 \\ 1\leq |m'|\ll q't^{1+\epsilon}/N\end{subarray}}{\sum\sum} \frac{1}{aq^{1-i\tau}} \frac{1}{aq'^{1+i\tau'}} I_1(q,m,\tau)\overline{I_1(q',m',\tau')} \textbf{T} d\tau d\tau'
\end{split}
\end{equation*}

where
\begin{equation*}
\textbf{T} = \sum_n n^{-1+\frac{-i\tau+i\tau'}{2}} U\left(\frac{n}{L}\right) e\left(\frac{n(a'q-aq')}{qq'}\right)
\end{equation*}

Analyzing $\textbf{T}$: Breaking the sum modulo $qq'$,
\begin{equation*}
\textbf{T} = \sum_{\beta (qq')} e\left(\frac{\beta(a'q-aq')}{qq'}\right) \sum_{l\in\Z}(\beta+qq'l)^{-1+\frac{-i\tau+i\tau'}{2}} U\left(\frac{\beta+qq'l}{L}\right)
\end{equation*}
applying Poisson summation to the $l-$sum,
\begin{equation*}
\textbf{T} = \sum_{\beta (qq')} e\left(\frac{\beta(a'q-aq')}{qq'}\right) \sum_{n\in\Z}\int_{\R} (\beta+qq'y)^{-1+\frac{-i\tau+i\tau'}{2}} U\left(\frac{\beta+qq'y}{L}\right)e(-ny)dy
\end{equation*}
and changing variables $w=(\beta + qq'y)/L$,
\begin{equation*}
\begin{split}
\textbf{T} &= \frac{1}{qq'}\sum_{\beta (qq')} e\left(\frac{\beta(a'q-aq')}{qq'}\right) \sum_{n\in\Z} e\left(\frac{n\beta}{qq'}\right)L^{(-i\tau+i\tau')/2}U^{\dagger}\left(\frac{nL}{qq'}, -\frac{i\tau}{2} + \frac{i\tau'}{2}\right) \\ &= \frac{L^{-i\tau/2+i\tau'/2}}{qq'} \sum_{n\in \Z} \mathfrak{C} U^{\dagger}\left(\frac{nL}{qq'}, -\frac{i\tau}{2} + \frac{i\tau'}{2}\right)
\end{split}
\end{equation*}
with $\mathfrak{C}$ as before. Since $|\tau-\tau'|\ll (NK)^{1/2}t^{\epsilon}/C$, the bound on $U^{\dagger}$ gives arbitrary saving for $|n|\gg C(NK)^{1/2}t^{\epsilon}/L$. We therefore get

\begin{lem}
\begin{equation}
S_{1,J}(N,C,L) \ll \frac{K}{NC^4} \underset{\begin{subarray}{c} C<q\leq2C, (m,q)=1 \\ 1\leq |m|\ll qt^{1+\epsilon}/N\end{subarray}}{\sum\sum} \underset{\begin{subarray}{c} C<q'\leq2C, (m',q')=1 \\ 1\leq |m'|\ll q't^{1+\epsilon}/N\end{subarray}}{\sum\sum} \sum_{|n|\ll C(NK)^{1/2}t^{\epsilon}/L} |\mathfrak{C}| |\mathfrak{K}| + O(t^{-2015})
\end{equation}
where
\begin{equation}
\mathfrak{K} = \underset{\R^2}{\int\int}(NL)^{-i\tau/2+i\tau'/2}\gamma(1+i\tau)\overline{\gamma(1+i\tau')}\frac{1}{q^{-i\tau}q'^{i\tau}}W_J(\tau)W_J(\tau') \mathcal{I}_1(q,m,\tau)\overline{\mathcal{I}_1(q',m',\tau')} U^{\dagger}\left(\frac{nL}{qq'}, -\frac{i\tau}{2} + \frac{i\tau'}{2}\right) d\tau d\tau'
\end{equation}
\end{lem}

Using the expression for $\mathcal{I}_1(q,m,\tau)$ as given in lemma (\ref{I**}), we get the expression

\begin{equation}
\begin{split}
\mathfrak{K} = \frac{|c_4|^2}{K^2}\underset{\R^2}{\int\int} & \gamma(1+i\tau)\overline{\gamma(1+i\tau')}W_J(q,m.\tau)\overline{W_J(q',m',\tau')}\frac{(LN)^{-i\tau/2+i\tau'/2}}{q^{-i\tau}q'^{i\tau'}}\left(-\frac{(t+\tau/2)q}{2\pi eNm}\right)^{-i(t+\tau/2)} \\ &\left(-\frac{(t+\tau'/2)q'}{2\pi eNm'}\right)^{i(t+\tau'/2)}U^{\dagger}\left(\frac{nL}{qq'}, -\frac{i\tau}{2} + \frac{i\tau'}{2}\right) d\tau d\tau' 
\end{split}
\end{equation}
where
\begin{equation*}
W_J(q,m,\tau) = \frac{1}{(t+\tau/2)^{1/2}}W_J(\tau)\left(-\frac{(t+\tau/2)q}{2\pi eNm}\right)^{3/2} V\left(-\frac{(t+\tau/2)q}{2\pi Nm}\right) \int_0^1 V\left(\frac{\tau}{2K} - \frac{(t+\tau/2)x}{Kma}\right) dx
\end{equation*}
%

Since $u^{3/2}V(u)\ll 1$ and $\tau\ll J\ll t^{1-\epsilon}$, it follows that
\begin{equation}
\frac{\partial}{\partial\tau}W_J(q,m,\tau) \ll \frac{1}{t^{1/2}|\tau|}
\end{equation}
We also note that the $x$-integral inside the expression of $W_J(q,m,\tau)$ contributes a factor of the size of its length, which is $\ll Kma/(t+\tau)$. Since $m\ll Ct^{1+\varepsilon}/N$ and $\tau\ll t$, the contribution is $\ll KCQt^{\varepsilon}/N$. Therefore $W_{J}(q,m,\tau)\ll K^{1/2}C/t^{1/2}N^{1/2}$.

We analyze the integral $\mathfrak{K}$ in two cases, when $n=0$ and when $n\neq0$. For $n=0$, the expression for $\mathfrak{C}$ gives $q=q'$, and the bound on $U^{\dagger}$ gives us arbitrary saving for $|\tau-\tau'|\gg t^{\epsilon}$. In this case,
\begin{equation*}
\mathfrak{K}\ll \frac{|c_4|^2}{K^2}\underset{|\tau|\ll (NK)^{1/2}/C}{\int}\gamma(1+i\tau)|^2W_J(q,m,\tau)\int_{|\tau'-\tau|\ll t^{\varepsilon}}W_J(q,m',\tau')d\tau' d\tau \ll \frac{t^{\epsilon}C}{K^{1/2}N^{1/2}t} =: B^*(C,0)
\end{equation*}
When $n\neq0$, 
\begin{equation}
\begin{split}
U^{\dagger}\left(\frac{nL}{qq'}, 1-\frac{i\tau}{2} + \frac{i\tau'}{2}\right) = \frac{c_5}{(\tau-\tau')^{1/2}} &U\left(\frac{(\tau-\tau')qq'}{4\pi nL}\right)\left(\frac{(\tau-\tau')qq'}{4\pi enL}\right)^{-i\tau/2+i\tau'/2} \\&+ O\left(\min\left\lbrace\frac{1}{|\tau-\tau'|^{3/2}}, \frac{C^3}{(|n|L)^{3/2}}\right\rbrace\right)
\end{split}
\end{equation}
for some absolute constant $c_5$.

Contribution of the error term towards $\mathfrak{K}$ is of the order of

\begin{equation*}
\frac{t^{\epsilon}}{K^2}\underset{[J,4J/3]^2}{\int\int} \frac{1}{t} \min\left\lbrace\frac{1}{|\tau-\tau'|^{3/2}}, \frac{C^3}{(|n|L)^{3/2}}\right\rbrace d\tau d\tau'
\end{equation*}
When the second term is smaller,
\begin{equation}
\frac{t^{\epsilon}}{K^2}\underset{\begin{subarray}{c} [J,4J/3]^2 \\ |\tau-\tau'|\ll |n|L/C^2\end{subarray}}{\int\int} \frac{1}{t} \frac{C^3}{(|n|L)^{3/2}} d\tau d\tau' \ll \frac{1}{K^{3/2}t}\frac{N^{1/2}}{(|n|L)^{1/2}}t^{\epsilon}
\end{equation}
When the first term is smaller,
\begin{equation}
\begin{split}
\frac{t^{\epsilon}}{K^2}\underset{\begin{subarray}{c} [J,4J/3]^2 \\ |\tau-\tau'|\gg |n|L/C^2\end{subarray}}{\int\int} \frac{1}{t} \frac{1}{|\tau-\tau'|^{3/2}} d\tau d\tau' & \ll \frac{t^{\epsilon}}{K^2t} \frac{C}{(|n|L)^{1/2}} \underset{[J,4J/3]^2}{\int\int}\frac{1}{|\tau-\tau'|^{1-\epsilon}}d\tau d\tau' \\ & \ll \frac{1}{K^{3/2}t}\frac{N^{1/2}}{(|n|L)^{1/2}}t^{\epsilon}
\end{split}
\end{equation}
The error contribution (for $n\neq0$) is
\begin{equation*}
B^*(C,n) = \frac{1}{K^{3/2}t}\frac{N^{1/2}}{(|n|L)^{1/2}}t^{\epsilon}
\end{equation*}


We finally analyze the main term. Striling's formula is
\begin{equation*}
\Gamma(\sigma+i\tau) = \sqrt{2\pi}(i\tau)^{\sigma-1/2}e^{-\pi|\tau|/2}\left(\frac{|\tau|}{e}\right)^{i\tau}\left\lbrace 1+ O\left(\frac{1}{|\tau|}\right) \right\rbrace
\end{equation*}
as $|\tau|\rightarrow\infty$. That gives
\begin{equation}
\gamma(1+i\tau) = \left(\frac{|\tau|}{4\pi e}\right)^{i\tau} \Phi(\tau), \quad \textit{where } \Phi'(\tau)\ll \frac{1}{|\tau|}
\end{equation}

By Fourier inversion, we write
\begin{equation*}
\left(\frac{4\pi nL}{(\tau-\tau')qq'}\right)^{1/2} U\left(\frac{(\tau-\tau')qq'}{4\pi nL}\right) = \int_{\R} U^{\dagger}(r,1/2)e\left(\frac{(\tau-\tau')qq'}{4\pi nL}r\right) dr
\end{equation*}
We conclude that for some constant $c_6$ (depending on the sign of $n$)
\begin{equation}\label{Kfinal}
\mathfrak{K} = \frac{c_6}{K^2}\left(\frac{qq'}{|n|L}\right)^{1/2}\int_{\R}U^{\dagger}(r,1/2)\underset{\R^2}{\int\int} g(\tau,\tau')e(f(\tau,\tau'))d\tau d\tau' dr + O(B^*(C,n))
\end{equation}
where
\begin{equation*}
\begin{split}
2\pi f(\tau,\tau') = & \tau\log\left(\frac{\tau}{4\pi e}\right) - \tau'\log\left(\frac{\tau'}{4\pi e}\right) - \frac{(\tau-\tau')}{2}\log(LN) + \tau\log q - \tau'\log q' \\ & -(t+\tau/2)\log\left(-\frac{(t+\tau/2)q}{2\pi eNm}\right) + (t+\tau'/2)\log\left(-\frac{(t+\tau'/2)q'}{2\pi eNm'}\right) \\ &\frac{(\tau-\tau')}{2}\log\left(\frac{(\tau-\tau')qq'}{4\pi enL}\right) + \frac{(\tau-\tau')qq'}{2nL}r
\end{split}
\end{equation*}
and
\begin{equation*}
g(\tau,\tau') = \Phi(\tau)\overline{\Phi(\tau')}W_J(q,m,\tau)W_J(q',m',\tau')
\end{equation*}

We intend to use the second derivative bound as given in Lemma \ref{2nd.der.lem}. For that, we need the following
\begin{equation*}
2\pi\frac{\partial^2}{\partial\tau^2}f(\tau,\tau') = \frac{1}{4}\left(\frac{4}{\tau} - \frac{1}{(t+\tau/2)} + \frac{2}{(\tau'-\tau)}\right), \quad 2\pi\frac{\partial^2}{\partial\tau'^2}f(\tau,\tau') = \frac{1}{4}\left(\frac{-4}{\tau'} + \frac{1}{(t+\tau'/2)} + \frac{2}{(\tau'-\tau)}\right)
\end{equation*}
and
\begin{equation*}
2\pi\frac{\partial^2}{\partial\tau'\partial\tau}f(\tau,\tau') = \frac{-1}{4}\left(\frac{2}{\tau'-\tau}\right)
\end{equation*}
Also, by explicit computation,
\begin{equation*}
4\pi^2\left[ \frac{\partial^2}{\partial\tau^2}f(\tau,\tau')\frac{\partial^2}{\partial\tau'^2}f(\tau,\tau') - \left(\frac{\partial^2}{\partial\tau'\partial\tau}f(\tau,\tau')\right)^2 \right] = -\frac{1}{2\tau\tau'} + O\left(\frac{1}{tJ}\right)
\end{equation*}
for $\tau,\tau'$ such that $g(\tau,\tau')\neq0$. So the conditions of lemma 4 of Munshi \cite{munshi} hold with $r_1=r_2=1/J^{1/2}$. To calculate the total variation of $g(\tau,\tau')$, recall that $\Phi'(\tau)\ll |\tau|^{-1}$ and $W'_J(q,m,\tau)\ll t^{-1/2}|\tau|^{-1}$. So $var(g)\ll t^{-1+\epsilon}$. So the double integral in (\ref{Kfinal}) over $\tau,\tau'$ is bounded by $O(Jt^{-1+\epsilon})$. Integrating trivially over $r$ using the rapid decay of the Fourier transform, we get that total contribution of the leading term in (\ref{Kfinal}) towards $\mathfrak{K}$ is bounded by
\begin{equation*}
O\left( \frac{1}{K^2} \frac{C}{(|n|L)^{1/2}} \frac{(NK)^{1/2}}{C}t^{-1+\epsilon}\right) = O(B^*(C,n))
\end{equation*}

Putting everything together, we get the final bound
\begin{equation*}
\begin{split}
S_{1,J}(N,C,L) & \ll \frac{t^{\epsilon}K}{NC^2}\bigg[ \underbrace{\underset{\begin{subarray}{c} C<q\leq 2C, (m,q)=1 \\ 1\leq m\ll qt^{1+\epsilon}/N\end{subarray}}{\sum\sum}}_{n=0} \left(\frac{t}{N}\right)B^*(C,0) +\underset{n\neq0}{\sum_{|n|\ll \frac{C(NK)^{1/2}t^{\epsilon}}{L}}}\sum_{C<q\leq 2C}\sum_{C<q'\leq 2C}\left(\frac{t}{N}\right)^2 B^*(C,n) \bigg] \\ & = \frac{t^{\epsilon}K}{NC^2} \bigg[\frac{C^3t}{N^{5/2}K^{1/2}} + \frac{C^{1/2}(NK)^{1/4}}{L}\frac{C^2t}{N^{3/2}K^{3/2}} \bigg]
\end{split}
\end{equation*}

That gives
\begin{equation*}
\begin{split}
S_{1,J}(N,C) & \leq t^{\epsilon}N^{1/2}K \underset{dyadic}{\sum_{1\leq L\ll Kt^{\epsilon}}} L^{1/2} \frac{K^{1/2}}{N^{1/2}C} \left[\frac{C^{3/2}t^{1/2}}{N^{5/4}K^{1/4}} + \frac{C^{1/4}(NK)^{1/8}}{L^{1/2}}\frac{Ct^{1/2}}{N^{3/4}K^{3/4}} \right] \\ &\ll t^{\epsilon}K^{3/2}\left( \frac{K^{1/4}C^{1/2}t^{1/2}}{N^{5/4}} + \frac{C^{1/4}t^{1/2}}{(NK)^{5/8}} \right)
\end{split}
\end{equation*}
Multiplying by $N^{1/2}/K$ and summing over the dyadic range $C\ll Q$, we get
\begin{equation}\label{S_1}
\frac{S_1(N)}{N^{1/2}} \ll t^{1/2+\varepsilon} \left( \frac{K^{1/2}}{N^{1/2}} + \frac{1}{K^{1/4}} \right)
\end{equation}
Finally, from equations (\ref{S_2}) and (\ref{S_1}), it follows that for $N\ll t^{1+\epsilon}$ and $K\gg N^{1/2}$,
\begin{equation*}
\frac{S(N)}{N^{1/2}}\ll t^{1/2+\varepsilon}\left(\frac{K^{1/2}}{N^{1/2}}+\frac{N^{1/4}}{K^{3/4}}+ \frac{K^{1/2}}{N^{1/2}} + \frac{1}{K^{1/4}} \right).
\end{equation*}
The optimal choice for $K$ occurs at $K=N^{2/3}$ and we get Proposition \ref{mainprop}.

\subsection*{Acknowledgments} I would like to thank Prof. Ritabrata Munshi for suggesting the problem, and Prof. Roman Holowinsky for many insightful discussions and encouragement.

\bibliographystyle{amsplain}
\bibliography{ref}

\providecommand{\bysame}{\leavevmode\hbox to3em{\hrulefill}\thinspace}
\providecommand{\MR}{\relax\ifhmode\unskip\space\fi MR }
\providecommand{\MRhref}[2]{%
  \href{http://www.ams.org/mathscinet-getitem?mr=#1}{#2}
}
\providecommand{\href}[2]{#2}
\begin{thebibliography}{10}

\bibitem{bky}
Valentin Blomer, Rizwanur Khan, and Matthew Young, \emph{Distribution of mass
  of holomorphic cusp forms}, Duke Math. J. \textbf{162} (2013), no.~14,
  2609--2644. \MR{3127809}

\bibitem{good82}
Anton Good, \emph{The square mean of {D}irichlet series associated with cusp
  forms}, Mathematika \textbf{29} (1982), no.~2, 278--295 (1983). \MR{696884}

\bibitem{huxley}
M.~N. Huxley, \emph{On stationary phase integrals}, Glasgow Math. J.
  \textbf{36} (1994), no.~3, 355--362. \MR{1295511}

\bibitem{ju97}
Matti Jutila, \emph{Mean values of {D}irichlet series via {L}aplace
  transforms}, Analytic number theory ({K}yoto, 1996), London Math. Soc.
  Lecture Note Ser., vol. 247, Cambridge Univ. Press, Cambridge, 1997,
  pp.~169--207. \MR{1694992}

\bibitem{kmv}
E.~Kowalski, P.~Michel, and J.~VanderKam, \emph{{Rankin-Selberg $L$-functions
  in the level aspect}}, Duke Math. J. \textbf{114} (2002), no.~1, 123--191.

\bibitem{li11}
Xiaoqing Li, \emph{Bounds for {${\rm GL}(3)\times {\rm GL}(2)$} {$L$}-functions
  and {${\rm GL}(3)$} {$L$}-functions}, Ann. of Math. (2) \textbf{173} (2011),
  no.~1, 301--336. \MR{2753605}

\bibitem{mive10}
Philippe Michel and Akshay Venkatesh, \emph{The subconvexity problem for {${\rm
  GL}_2$}}, Publ. Math. Inst. Hautes \'Etudes Sci. (2010), no.~111, 171--271.
  \MR{2653249 (2012c:11111)}

\bibitem{munshi}
Ritabrata Munshi, \emph{The circle method and bounds for
  {$L$}-functions--{III}: $t-$aspect subconvexity for {$GL(3) L$}-functions},
  J. Amer. Math. Soc. \textbf{28} (2015), 913--938. \MR{3369905}

\bibitem{mu15}
\bysame, \emph{The circle method and bounds for {$L$}-functions, {II}:
  {S}ubconvexity for twists of {${\rm GL}(3)$} {$L$}-functions}, Amer. J. Math.
  \textbf{137} (2015), no.~3, 791--812. \MR{3357122}

\bibitem{singh17}
S.~K. {Singh}, \emph{{$\backslash$lowercase$\{$t$\}$-$\backslash$lowercase
  $\{$aspect subconvexity bound for$\}$ $GL(2)$
  L-$\backslash$lowercase$\{$functions $\}$}, arxiv:1706.04977}, June 2017.

\bibitem{srini65}
B.~R. Srinivasan, \emph{The lattice point problem of many dimensional
  hyperboloids. {III}}, Math. Ann. \textbf{160} (1965), 280--311. \MR{0181614}

\end{thebibliography}





\address{Department of Mathematics, The Ohio State University, 100 Math Tower,
231 West 18th Avenue, Columbus, OH 43210-1174}\\
Email address: \email{aggarwal.78@osu.edu}

\end{document}